\documentclass[preprint,11pt]{elsarticle}
\biboptions{numbers,sort&compress}
\usepackage{amsfonts,psfrag,amsmath,mathrsfs,amsthm,color,empheq}
\usepackage{color,amssymb,bm} 
\usepackage{epsfig,soul}
\usepackage[most]{tcolorbox}

\newtcolorbox{tcbdoublebox}[1][]{%
	enhanced jigsaw,
	sharp corners,
	colback=white,
	fontupper={\setlength{\parindent}{20pt}},
	#1
}

\usepackage{hyperref}
\hypersetup{
	colorlinks   = true, 
	urlcolor     = blue, 
	linkcolor    = blue, 
	citecolor   = red 
}
\usepackage{multirow}

 \renewenvironment{thebibliography}[1]{%
 	\begin{oldthebibliography}{#1}%
 		\setlength{\parskip}{0.1ex}%
 		\setlength{\itemsep}{0.1ex}%
 	}%
 	{%
 	\end{oldthebibliography}%
 }

\usepackage{mathptmx} 
\usepackage{amsmath,mathtools} 
\usepackage{amsthm,amssymb,bbm}

 \usepackage{graphicx,float}
 \usepackage[export]{adjustbox}
 \usepackage{setspace}
 \usepackage{newtxmath}
 \usepackage{mathtools} 
 
\usepackage[T1]{fontenc}
\usepackage{tgschola}

\usepackage{setspace}
\AtBeginDocument{%
	\addtolength\abovedisplayskip{-0.2\baselineskip}%
	\addtolength\belowdisplayskip{-0.2\baselineskip}%
}
 \usepackage[margin=1in]{geometry}

\newcommand{\IC}{{\mathbb{C}}}

\newcommand{\CZ}{{\mathcal{Z}}}

\newcommand{\ID}{{\mathbb{D}}}

\newcommand{\diag}{{\rm{diag}}}

\newcommand{\bmt}{\left[ \begin{array}{ccccccccc}}
\newcommand{\emt}{\end{array}\right]}
\newcommand{\bean}{\begin{eqnarray*}}
\newcommand{\eean}{\end{eqnarray*}}
\newcommand{\bea}{\begin{eqnarray}}
\newcommand{\eea}{\end{eqnarray}}
\newcommand{\eq}{\begin{equation}\begin{array}{lllllllll}}
\newcommand{\ee}{\end{array}\end{equation}}
\newcommand{\eqn}{\begin{equation*}\begin{array}{lllllllll}}
\newcommand{\een}{\end{array}\end{equation*}}

\allowdisplaybreaks
\usepackage[symbol]{footmisc} 

\graphicspath{{./Figures/}}

\journal{ERA}

\begin{document}

	\begin{frontmatter}
	
	\title{A regularized eigenmatrix method for unstructured sparse recovery
	}
	
	 	\author{Koung Hee Leem} 	
	\author{Jun Liu\corref{mycorrespondingauthor}}
		\cortext[mycorrespondingauthor]{Corresponding author.}
	\ead{juliu@siue.edu}  	
	\author{George Pelekanos} 	

	\address{Department of Mathematics and Statistics, Southern Illinois University Edwardsville, Edwardsville, IL 62026, USA.\vspace{-1cm}}

	\begin{abstract}
	The recently developed data-driven eigenmatrix method shows very promising reconstruction accuracy in sparse recovery for a wide range of kernel functions and random sample locations.
	However, its current implementation can lead to numerical instability if the threshold tolerance is not appropriately chosen. To incorporate regularization techniques, we propose to regularize the eigenmatrix method by replacing the computation of an ill-conditioned pseudo-inverse by the solution of an ill-conditioned least square system, which can be efficiently treated by Tikhonov regularization. Extensive numerical examples confirmed the improved effectiveness of our proposed method, especially when the noise levels are relatively high.
	\end{abstract}
	
	
\end{frontmatter}

\section{Background} 
Let $X$ be the parameter space and $S$ be the sampling space. Assume $g(s,x)$ is a given kernel function on $S\times X$ that is analytic in $x$. Suppose the unknown sparse signal $f$ is given by
\eq
f(x)=\sum_{k=1}^{n_x} w_k \delta(x-x_k),
\ee
with $\delta$ being the Dirac delta function, and  $n_x$ spikes with distinct locations $\{x_k\}$ and weights $\{w_k\}$. 
The observable for any given sampling point $s\in S$ is given by the following summation
\eq \label{u_obs}
u(s):=\int_{X} g(s,x)f(x)\ dx=\sum_{k=1}^{n_x} w_k g(s,x_k).
\ee 
Let $\{s_j\}$ be a chosen set of $n_s$ (unstructured) sample locations in $S$ and $u_j=u(s_j)$ be the unknown exact values of observations. In practice, we can only obtain noisy observations, which are assumed to have the following multiplicative form (with an unknown noise magnitude $\sigma>0$)
\eq
\widetilde{u}_j=u_j(1+\sigma \CZ_j)
\ee 
with $\CZ_j$ being  independently identically distributed (i.i.d.) standard Gaussian random variables with zero mean and unit variance. Our task is to recover the unknown spike locations $\bm x:=[x_1;x_2;\cdots;x_{n_x}]$ and weights $\bm w:=[w_1;w_2;\cdots;w_{n_x}]$ from the observation $\{\widetilde{u}_j\}$. Obviously, this leads to a highly nonlinear inverse problem that is difficult to treat numerically. The standard nonlinear least square formulation will lead to a nonconvex unconstrained optimization problem that can be better solved with a good initial guess estimated by the proposed methods.

Depending on the definition of kernel function $g$, the sparse recovery problem in the above general form (\ref{u_obs}) covers a list of well-known sparse recovery problems,
such as rational approximation \cite{berljafa2017rkfit}, spectral function estimation \cite{ying2022analytic,ying2022pole}, Fourier inversion \cite{potts2010parameter,potts2013parameter}, Laplace inversion \cite{weeks1966numerical,davies1979numerical,peter2013generalized,cohen2007numerical}, and sparse deconvolution, for which many specially designed numerical algorithms \cite{becker2011nesta,marques2018review} were established with sounding theoretical support in the past few decades; see references in \cite{ying2024eigenmatrix}. Nevertheless, these tailored algorithms rely heavily on the underlying structure of each problem, which are not directly applicable to general  kernel function with unstructured sampling grid. The developed data-driven eigenmatrix method in \cite{ying2024eigenmatrix} does not assume any structures in the kernel function and sampling grid and hence it has a wider applicability that specialized or structured sparse recovery algorithms. Nevertheless, it requires the computation of the psedudo-inverse of a highly ill-conditioned rectangle matrix, which can lead to numerical instability when the threshold tolerance does not match with the underlying noise levels in the measurement data. 
Our major contribution is to propose a regularized eigenmatrix method that can handle noisy measurement data through modern Tikhonov regularization techniques, which demonstrates significantly improved recovery accuracy in tested numerical examples with high noise levels.

This paper is organized as follows. In
Section 2 we briefly review the original eigenmatrix method and point out its drawbacks. In Section 3 we introduce a new regularized eigenmatrix method based on Tikhonov regularization techniques.
A few numerical examples are presented in Section 4. Finally, some remarks are concluded in Section 5.
\section{Review of the eigenmatrix method}
Inspired by the shifting operator defined in the Prony’s method and the ESPRIT algorithm \cite{roy1989esprit},
the recently developed eigenmatrix method \cite{ying2024eigenmatrix} for unstructured sparse recovery problems shows very appealing reconstruction accuracy for different kernels and unstructured sampling locations.
Its key idea is to find an $n_s$-by-$n_s$ eigenmatrix $M$ such that for all $x\in X$ there approximately holds \textit{eigensytem}
\eq \label{Meig}
M \bm g(x)\approx x \bm g(x),
\ee
where  
$
\bm g(x)=[g(s_j,x)]_{1\le j\le n_s}
$ is an $n_s$-by-1 vector of functions on $X$. In numerical implementations, we can enforce this approximate relation over a set of collocation nodes $\{a_t\}_{t=1}^{n_a}$ selected in $X$. 
More specifically, if $X$ is the unit disk $\ID$ on complex plane, one can select a uniform grid of collocation nodes on the boundary of unit disk, which can be justified by invoking the exponential convergence of trapezoidal rule and the application of Cauchy integral theorem for analytic functions. 
If $X$ is the real interval $[-1,1]$, one can choose a Chebyshev grid of collocation nodes on $[-1,1]$, which can be explained by the Chebyshev quadrature for analytic functions. A general connected domain $X$ can be treated by introducing a smooth one-to-one map between $X$ and $\ID$ or $[-1,1]$.
At this point, there is no error estimates on the accuracy of the approximation (\ref{Meig}) in various settings.

Following the notation and methodology  introduced in \cite{ying2024eigenmatrix}, the original eigenmatrix method based on the ESPRIT method mainly consists of the following 4 major steps (not including the postprocessing step for simplicity):
\begin{tcbdoublebox}[title={The original eigenmatrix method \cite{ying2024eigenmatrix}}]
\begin{enumerate}
	\item[(1)] Compute the matrix $G=[g(s_j,a_t)]\in\IC^{n_s\times n_a}$ based on the $n_s$ sampling locations $\{s_j\}_{j=1}^{n_s}$ and $n_a$ collocation nodes $\{a_t\}_{t=1}^{n_a}$. Normalize $G$ column-wisely to get $\widehat{G}$.
	\item[(2)] Compute the $n_s\times n_s$ eigenmatrix $M=\widehat{G}\Lambda \widehat{G}^\dagger$,
	where $\Lambda=\diag(a_t)$ and $\widehat{G}^\dagger$ is the pseudo-inverse of $\widehat{G}$ by thresholding singular values smaller than a given tolerance $tol$.
	\item[(3)] Given the vector of noisy observations $\widetilde{\bm u}$,  
		 choose $l>n_x$ and then compute rank-$n_x$ truncated SVD of the following matrix
			\eq \label{Ap2}
		 A:=\bmt \widetilde{\bm u},& M \widetilde{\bm u},&\cdots,&M^{l}\widetilde{\bm u} \emt  =U S V^*.
		\ee
		Define $V_+^*$ and $V_-^*$ be the sub-matrix of $V^*$ by deleting the first column and the last column, respectively. The $n_x$ eigenvalues $\{\widetilde{x}_k\}$ of the matrix $V_+^*(V_-^*)^\dagger $ yield the estimated spike locations. Here, we expect $(V_-^*)$ to be well-conditioned.
	 
	\item[(4)] With computed $\{\widetilde{x}_k\}$, the weights $\widetilde{\bm w}=[\widetilde{w}_1;\widetilde{w}_2;\cdots;\widetilde{w}_{n_x}]$ can be estimated via  a least square problem defined by $\widetilde{G} \widetilde{\bm w}=\widetilde{\bm u}$,
	where $\widetilde{G}=[g(s_j,\widetilde{x}_k)]$ is of size $n_s\times n_x$.
\end{enumerate}
\end{tcbdoublebox}
As a data-driven approach, it involves the key procedure of (approximately) finding the pseudo-inverse $\widehat{G}^\dagger$ of a highly ill-conditioned rectangular matrix $\widehat{G}$, which was not carefully treated from the perspective of regularization. To alleviate the issue of large condition number, in \cite{ying2024eigenmatrix} the author  suggested to choose (a small) $n_a=32$ such that $\widehat{G}$ is of full column rank and its condition number  is bounded below by $10^7$.
Moreover, the selected thresholding tolerance $tol$ was such that $\|M\|$ is bounded by a small constant such as $3$. In their implementations  \footnote{\url{https://github.com/lexingying/EigenMatrix}} however,  the author used $tol=10^{-4}\|\widehat G\|_{F}$ or  $tol=10^{-8}\|\widehat G\|_{F}$ as the thresholding tolerance in different examples. Hence, the strategy of selecting a small $n_a$ and $tol$, essentially points to some heuristic regularization treatment, which cannot take into account the actual noise level in the measurements. Therefore, the current version of the eigenmatrix method is less robust in handling a wide range of unknown noise levels.

\section{A regularized eigenmatrix method}

To make use of modern regularization techniques in the above eigenmatrix method, we need to avoid explicitly computing the ill-conditioned pseudo-inverse matrix $\widehat{G}^\dagger$.
In view of the matrix $A$ in (\ref{Ap2}), we only need the matrix-vector products $M^k \widetilde{\bm u} $ for $k\ge 1$, which implies that the explicit construction of matrix $M$ is unnecessary.
By the theory of pseudo-inverses, if $\widehat{G}$ has \textbf{linearly independent columns}, then there holds $\widehat{G}^\dagger \widehat{G}=I_{n_a}$, which leads to
\[M^k \widetilde{\bm u}=(\widehat{G}\Lambda \widehat{G}^\dagger)^k  \widetilde{\bm u}=\underbrace{\widehat{G}\Lambda \widehat{G}^\dagger \widehat{G}\Lambda \widehat{G}^\dagger \cdots \widehat{G}\Lambda \widehat{G}^\dagger}_{k\ \textnormal{times}} \widetilde{\bm u}=\widehat{G}\Lambda^k \widehat{G}^\dagger \widetilde{\bm u}.
\]
Let $\bm v=\widehat{G}^\dagger \widetilde{\bm u}$, then we can rewrite the matrix $A$ in  (\ref{Ap2}) in the form
	\eq \label{Ap22}
A  =\bmt \widetilde{\bm u},& \widehat{G}\Lambda \bm v,&\cdots,& \widehat{G}\Lambda^{l}\bm v\emt.
\ee
The vector $\bm v=\widehat{G}^\dagger \widetilde{\bm u}$ can then be obtained from the following ill-conditioned linear system 
 \eq  
\label{Gv=u} \widehat{G} \bm v=\widetilde{\bm u} 
\ee
since then we obtain
$
\label{Gv2} \bm v=I_{n_a} \bm v=\widehat{G}^\dagger (\widehat{G} \bm v)=\widehat{G}^\dagger \widetilde{\bm u}.
$
From the above it follows, that there is not  even a need to approximately compute the pseudo-inverse $\widehat{G}^\dagger$ or construct the matrix $M$.

In summary, we propose the following regularized eigenmatrix method without $M$:
\begin{tcbdoublebox}[title={Our proposed regularized eigenmatrix method}]
\begin{enumerate}
	\item[(1)] Unchanged.
	\item[(2)]  Solve the system (\ref{Gv=u}) for $\bm v$ by Tikhonov regularization method (see below).
	\item[(3)]  
		Construct $A$ using (\ref{Ap22}), and leave the remaining parts unchanged. 
	\item[(4)] Unchanged.
\end{enumerate}
\end{tcbdoublebox}
We reiterate here, that the significant improvement from the original eigenmatrix method is to avoid explicitly computing the eigenmatrix $M$ that requires the computation of the pseudo-inverse $\widehat{G}^\dagger$.
Moreover, the noisy observation $\widetilde{\bm u}$ will influence the computation of $\bm v$ through the Tikhonov regularization techniques, which do not rely on any manual adjustments of algorithmic parameters.

To demonstrate the role of  vector $\bm v$ in the above method,
we consider $\widetilde{\bm u}= \widehat{G} \bm v$,  and we let $\bm a^k=\mathrm{diag}(\Lambda^k):=[a_1^k;a_2^k;\cdots;a_{n_a}^k]$ be a column vector.  We then obtain the  factorization
\eq \label{ApImproved3}
A = \widehat{G}\bmt \bm v,&  \Lambda \bm v,&\cdots,&  \Lambda^{l}\bm v\emt 
=\widehat{G}\ \diag(\bm v)\ [\bm{1},\bm a,\cdots,\bm a^{l}],
\ee
which closely fits the desired factorization structure of the eigenmatrix method, that is
\eq \label{ApImproved4}
A \approx   \widetilde{G}\ \diag(\widetilde{\bm w})\ [\bm{1},\bm x,\cdots,\bm x^{l}].
\ee
Hence, the entries of the vector $\bm v$ act as the 'weights' for the corresponding collocation nodes $\{a_t\}$. This connection may be helpful to choose better sampling points $\{s_j\}$ and collocation nodes $\{a_t\}$. 
For instance, if  $n_a=n_x$ such that collocation nodes are identical to the  unknown spike locations, 
we would expect to obtain very accurate reconstruction.
\subsection{Tikhonov regularization for solving (\ref{Gv=u})}
We are now ready to employ modern regularization techniques \cite{engl1996regularization,hansen2010discrete,kirsch2011introduction,benning2018modern} to solve (\ref{Gv=u}). 
The standard Tikhonov regularization approximates the solution $\bm v$ of (\ref{Gv=u}) by ${\bm v}_\gamma$, which is given as the minimizer of the following Tikhonov regularized objective functional
\begin{equation} \label{Tikfun}
	{
		{\bm v}_\gamma:=\mathop{\arg\min}_{\bm v }\ \left(\| \widehat{G}  {\bm v }-\widetilde{\bm u}\|^2+\gamma^2 \|  \bm v\|^2\right),}
\end{equation}
where $\gamma>0$ is a regularization parameter to be determined.
The above Tikhonov minimization problem (\ref{Tikfun}) is mathematically equivalent to  solving the regularized normal equation
\begin{equation} \label{Tikregsys}
	(\widehat{G}^*\widehat{G} +{\gamma^2} I) {\bm v}_{\gamma}=\widehat{G}^* \widetilde{\bm u}.
\end{equation}   
There are many different a priori or a posteriori methods of choosing a good regularization parameter $\gamma>0$, such as the Morozov's discrepancy principle \cite{engl1987discrepancy} that requires the knowledge of noise level. In our numerical experiments, we will apply and compare the established IMPC \cite{bazan2009improved} and L-curve \cite{hansen1992analysis,engl1994using} technique\footnote{\url{http://www2.compute.dtu.dk/~pcha/Regutools/}} for estimating the regularization parameter $\gamma$, since both methods do not require a priori knowledge of the (unknown) noise level in the measured data $\widetilde{\bm u}$.
Both methods yield regularization parameters that are very close to each other, and hence deliver similar reconstruction accuracy. For large-scale problems, computationally more efficient iterative regularization techniques \cite{hansen2018air,gazzola2019ir} may also be used. 
Our major contribution is not to develop a new regularization method, but to reformulate the original eigenmatrix method such that the modern regularization techniques can be seamlessly employed.
\section{Numerical results}
In this section, we will numerically compare the original eigenmatrix method based on pseudo-inverse (denoted by pinv) and our proposed regularized eigenmatrix method based on IMPC and L-curve techniques. 
All simulations are implemented using MATLAB R2024a.
To better illustrate the influence of our proposed regularization techniques on the reconstruction accuracy,
we will only compare the recovered spike locations and weights based on the ESPRIT algorithm, without the extra postprocessing step of nonlinear optimization that may further improve the accuracy. 
To measure the reconstruction accuracy, we report the absolute difference of the spike locations and weights separately in Euclidean norm as the overall reconstruction errors, that is
\[
\textnormal{errors}=(\|\bm x-\widetilde{\bm x}\|_2,\|\bm w-\widetilde{\bm w}\|_2).
\]
We would expect the obtained errors to get smaller as the noise level $\delta$ is decreased. We will test the same examples and sampling locations as described in \cite{ying2024eigenmatrix}, except that the used noise levels are increased by 10 times to better demonstrate the robustness of our proposed regularization techniques. In particular, the spike weights $\{w_k\}$ are set to be one and the noisy observation $\widetilde{u}_j$ is constructed by adding different levels of Gaussian noise to the exact observation $u_j=u(s_j)$.  
The obtained errors are affected by both the used algorithms and the measurement noise, which may show some variance in numerical simulations. To minimize the influence of randomness, we compare all algorithms with the same random noise for a given noise level $\delta$.
\subsection{Example 1 (Rational approximation)}
\begin{figure}[htp!]
	\begin{center}
		\includegraphics[width=1\textwidth]{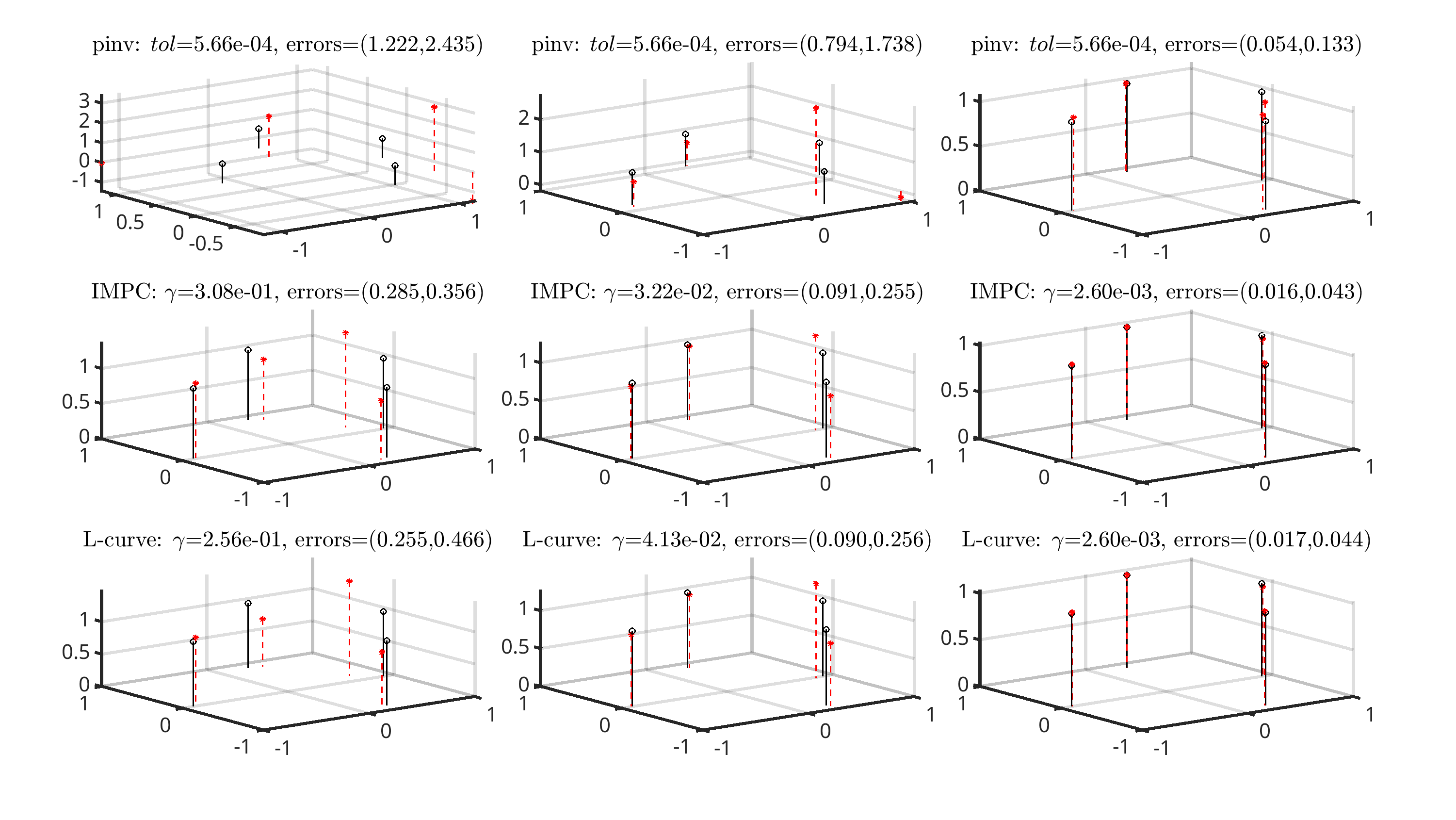} 
	\end{center}
\vspace{-30pt}
	\caption{Rational approximation (Columns from left to right: $\sigma=10^{-1},10^{-2},10^{-3}$). The exact spikes in solid line and the recovered spikes in dashed line. The errors measure the 2-norm difference in spike locations and weights, respectively.}	
	\vspace{-10pt}
	\label{FigEx1}
\end{figure}

In this problem, we have $X=\ID$, $g(s,x)=\frac{1}{s-x}$, and true spike locations $$\bm x=0.9 e^{2\pi i [0.2;0.5;0.8;1]}.$$ We generated $n_s=40$ random sampling points $\{s_j\}$ outside the unit disk, each with a modulus between 1.2 and 2.2. We then build the matrix $G=[g(s_j,a_t)]$ with $n_a=32$ uniformly spaced collocation nodes on the unit circle. Figure \ref{FigEx1} shows the reconstructed spike locations and weights in comparison with the exact ones by 3 different methods (from top to bottom: pinv, IMPC, and L-curve) with 3 different noise levels. 
Clearly, our regularized eigenmatrix methods (both IMPC and L-curve) deliver improved recovery (with smaller errors in each column), especially when the noise level gets higher. Both IMPC and L-curve technique yields comparable regularization parameter $\gamma$.
Notice the threshold tolerance $tol=10^{-4}\|\widehat{G}\|_F$ used in the original eigenmatrix method is independent of the noise level, which may cause degraded reconstruction accuracy if not appropriately chosen.

\subsection{Example 2 (Spectral function approximation)}
\begin{figure}[htp!]
	\begin{center}
		\includegraphics[width=1\textwidth]{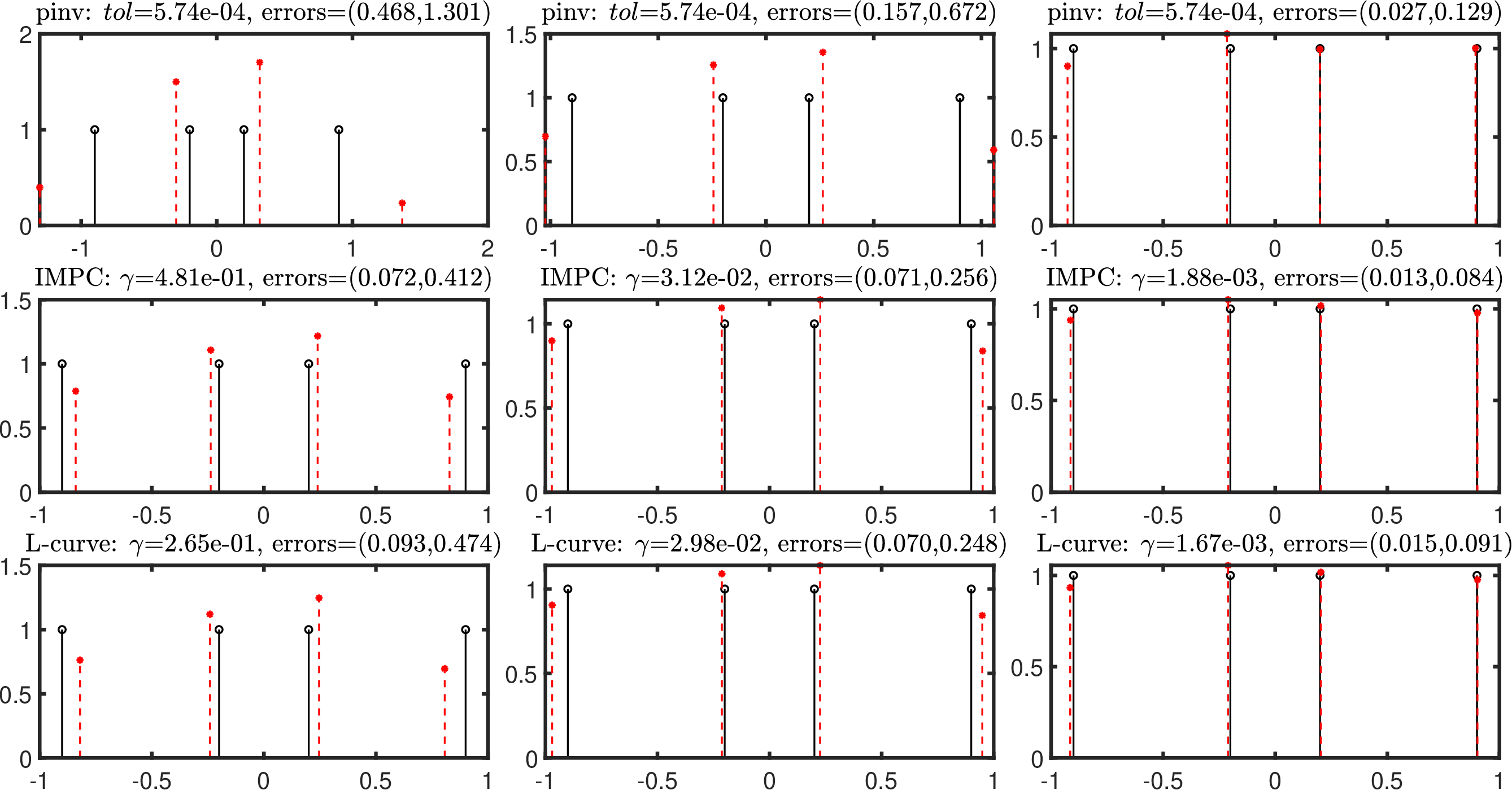} 
	\end{center}
	\caption{Spectral function approximation (Columns from left to right:$\sigma=10^{-1},10^{-2},10^{-3}$). The exact spikes in solid line and the recovered spikes in dashed line. The errors measure the 2-norm difference in spike locations and weights, respectively. }	
	\label{FigEx2}
\end{figure}

In this problem, we have $X=[-1,1]$, $g(s,x)=\frac{1}{s-x}$, and true spike locations $$\bm x=[-0.9; -0.2; 0.2; 0.9].$$ We use  $n_s=256$ uniformly distributed sampling points $s_j=\pm (2j-1)\pi i/\beta, j=1,2,\cdots,128$ from the Matsubara grid on the imaginary axis, and then build the matrix $G=[g(s_j,a_t)]$ with $n_a=32$ Chebyshev collocation nodes on $[-1,1]$. 
Figure \ref{FigEx2} shows the reconstructed spike locations and weights in comparison with the exact ones by 3 different methods (denoted by pinv, IMPC, and L-curve) with 3 different noise levels. Clearly, our regularized eigenmatrix method (both IMPC and L-curve) delivers more accurate recovery, especially when the noise levels are high.

\subsection{Example 3 (Fourier inversion)}
\begin{figure}[htp!]
	\begin{center}
		\includegraphics[width=1\textwidth]{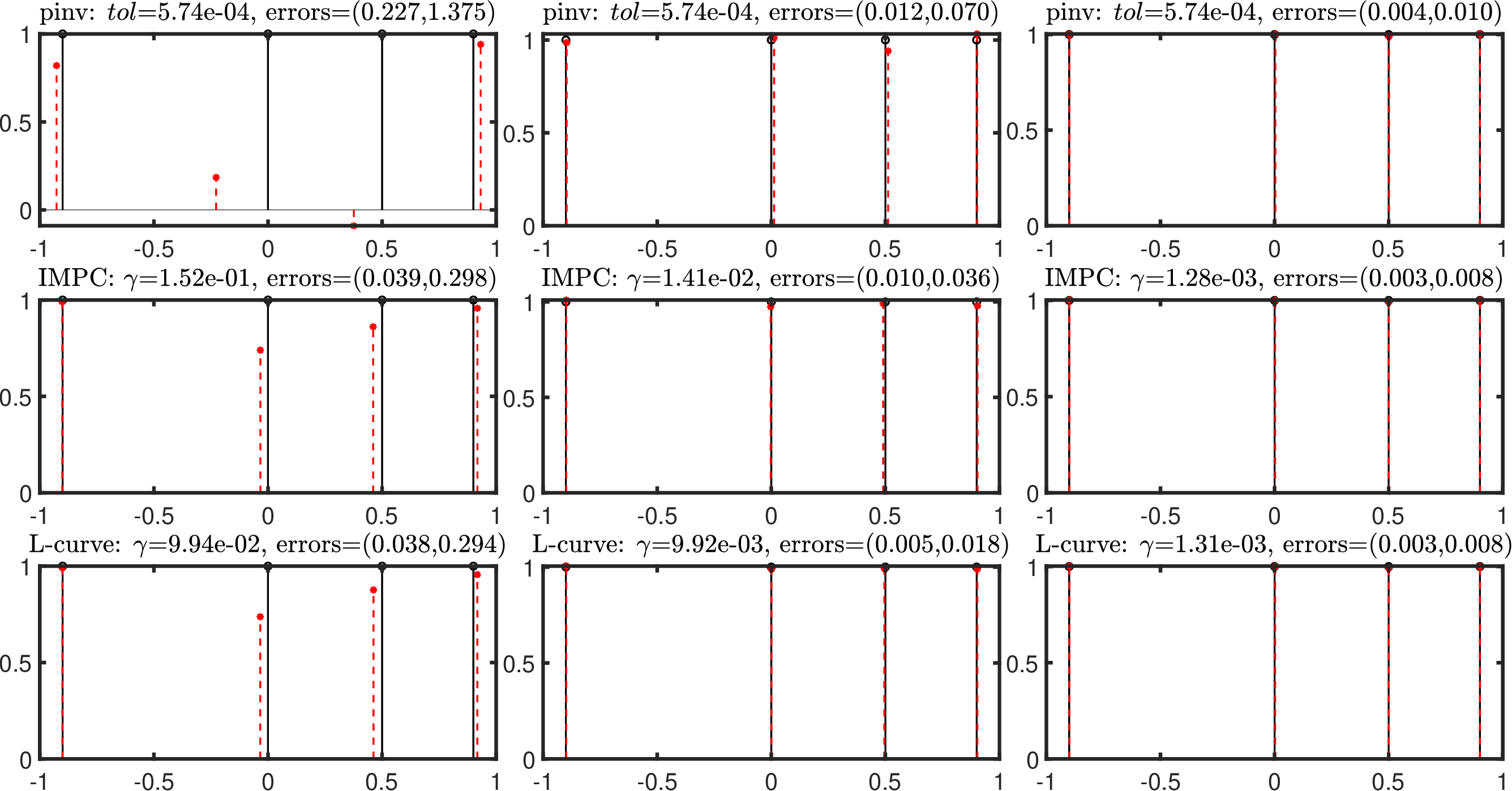} 
	\end{center}
	\caption{Fourier inversion (Columns from left to right:$\sigma=10^{-1},10^{-2},10^{-3}$). The exact spikes in solid line and the recovered spikes in dashed line. The errors measure the 2-norm difference in spike locations and weights, respectively. }	
	\label{FigEx3}
\end{figure}

In this problem, we have $X=[-1,1]$, $g(s,x)=e^{i\pi s x}$, and true spike locations $$\bm x=[-0.9;0;0.5;0.9].$$ We generated $n_s=128$ random sampling points $\{s_j\}$ in $[-5,5]$, and then build the matrix $G=[g(s_j,a_t)]$ with $n_a=32$ Chebyshev collocation nodes on $[-1,1]$. 
Figure \ref{FigEx3} presents the reconstructed spike locations and weights in comparison with the exact ones by 3 different methods (denoted by pinv, IMPC, and L-curve) with 3 different noise levels. Again, our regularized eigenmatrix method (both IMPC and L-curve) delivers more accurate recovery. It is worthwhile to point out that the original eigenmatrix method also works very well   
for this problem with small noise levels, likely due to the relatively smaller condition number $\mathrm{cond}(\widehat{G})\approx 10^{7}$.

\subsection{Example 4 (Laplace inversion)}
In this problem, we have $X=[0.1,2.1]$, $g(s,x)=xe^{-sx}$, and true spike locations $$\bm x=[0.2;1.1;1.6;2.0].$$ We generated $n_s=100$ random sampling points $\{s_j\}$ in $[0,10]$ and then build the matrix $G=[g(s_j,a_t)]$ with $n_a=32$ shifted Chebyshev collocation nodes on $[0.1,2.1]$. Notice here $\widehat{G}$ is not of full column rank with $\mathrm{rank}(\widehat{G})=17$ and a large condition number $\mathrm{cond}(\widehat{G})\approx 10^{17}$.  In \cite{ying2024eigenmatrix} the author only tested very low noise levels ($\sigma=10^{-5},10^{-6},10^{-7}$) with a very small threshold tolerance $tol=10^{-8}\|\widehat{G}\|_F$, which may conceal the essential difficulty of highly ill-conditioned $G$.  
Hence, we will test with higher noise levels ($\sigma=5\times 10^{-2},5\times 10^{-3},5\times 10^{-4}$),
where we found the moderate threshold tolerance $tol=10^{-4}\|\widehat{G}\|_F$ works better in treating higher noise levels.

\begin{figure}[htp!]
	\begin{center}
		\includegraphics[width=1\textwidth]{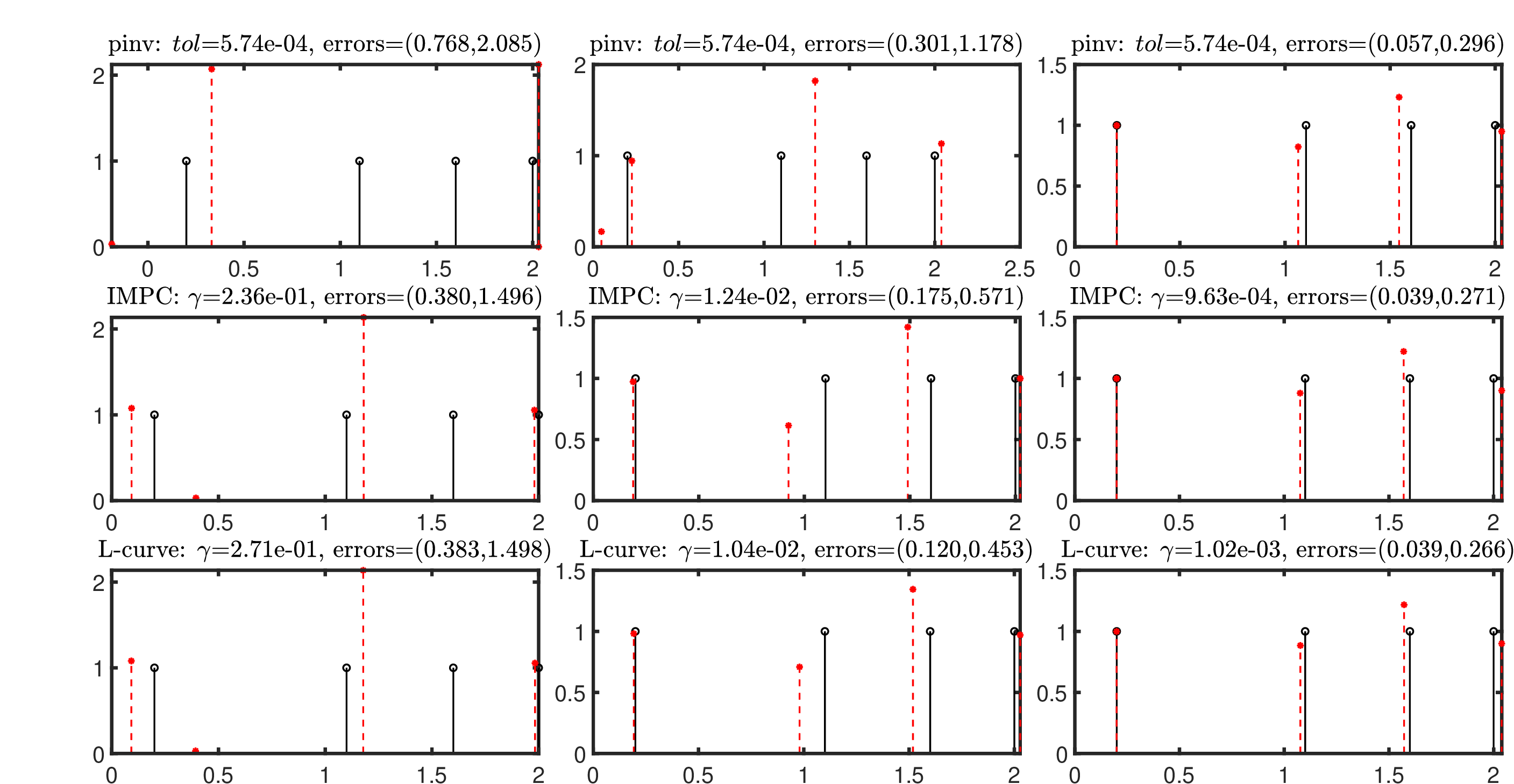} 
	\end{center}
	\vspace{-20pt}
	\caption{Laplace inversion (Columns from left to right: $2\sigma=10^{-1}, 10^{-2}, 10^{-3}$). The exact spikes in solid line and the recovered spikes in dashed line. The errors measure the 2-norm difference in spike locations and weights, respectively.}	
	\vspace{-10pt}
	\label{FigEx4}
\end{figure}
Figure \ref{FigEx4} shows the reconstructed spike locations and weights in comparison with the exact ones by 3 different methods (denoted by pinv, IMPC, and L-curve) with 3 different noise levels. 
Again, our Tikhonov regularized eigenmatrix method (both IMPC and L-curve) delivers more accurate recovery, which is expected since the most appropriate choice of an threshold tolerance $tol$ requires carefully tuning by hands.  

\subsection{Example 5 (Sparse deconvolution)}
\begin{figure}[htp!]
	\begin{center}
		\includegraphics[width=1\textwidth]{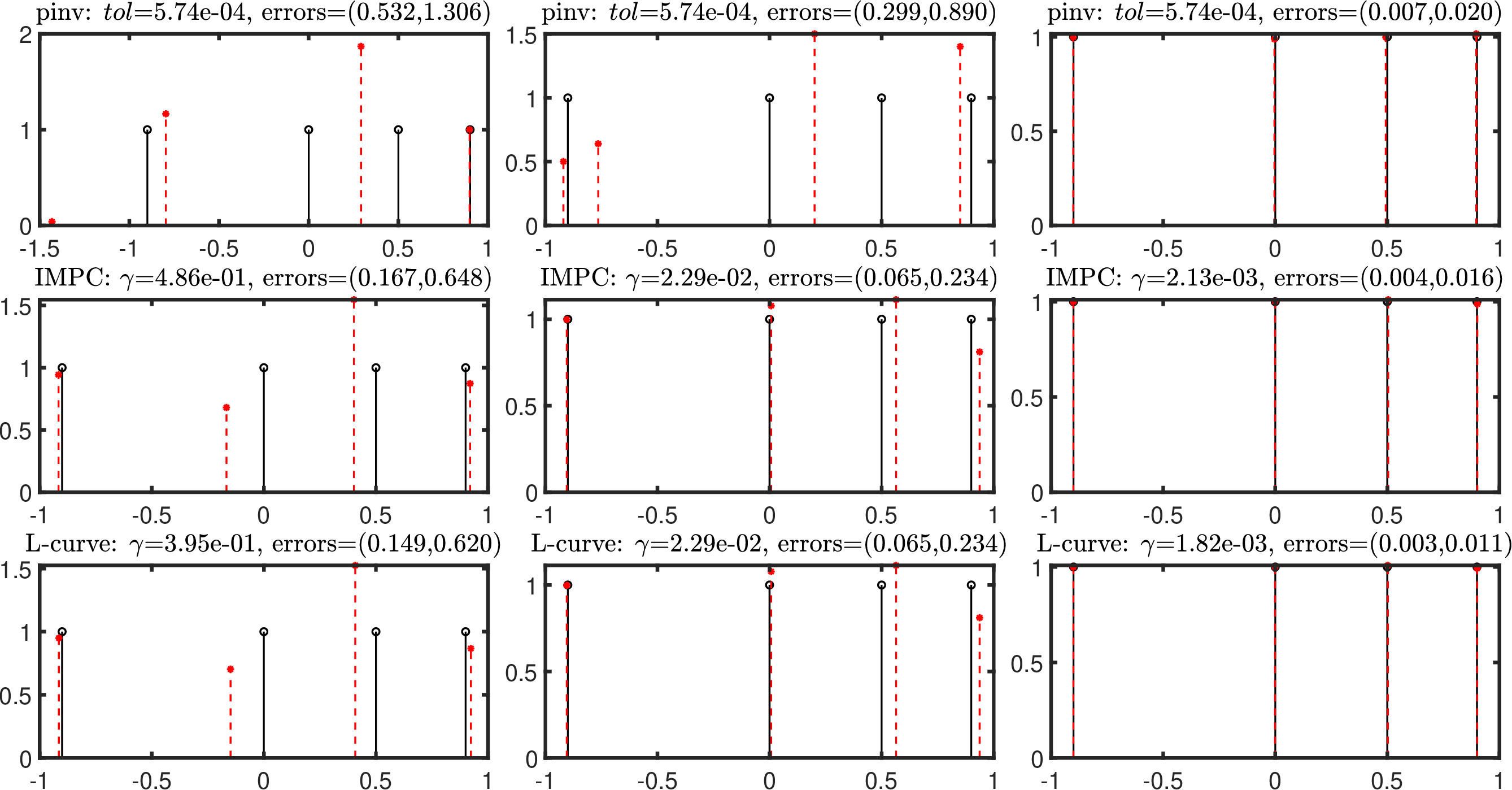} 
	\end{center}
	\caption{Sparse deconvolution (Columns from left to right:$\sigma=10^{-1},10^{-2},10^{-3}$). The exact spikes in solid line and the recovered spikes in dashed line. The errors measure the 2-norm difference in spike locations and weights, respectively.}	
	\label{FigEx5}
\end{figure}

In this problem, we have $X=[-1,1]$, $g(s,x)=\frac{1}{1+4(s-x)^2}$, and true spike locations $$\bm x=[-0.9;0;0.5;0.9].$$ We generated $n_s=128$ random sampling points $\{s_j\}$ in $[-5,5]$, and then build the matrix $G=[g(s_j,a_t)]$ with $n_a=32$ Chebyshev collocation nodes on $[-1,1]$. 
Figure \ref{FigEx5} displays the reconstructed spike locations and weights in comparison with the exact ones by 3 different methods (denoted by pinv, IMPC, and L-curve) with 3 different noise levels. Again, our Tikhonov regularized eigenmatrix method (both IMPC and L-curve) delivers more accurate recovery.

\section{Conclusion}
The original eigenmatrix method requires the computation of a pseudo-inverse matrix based on a  chosen threshold tolerance, which can not take into account the noise in data. Our proposed regularized eigenmatrix method addressed this shortcoming by incorporating modern regularization techniques, which provide improved recovery as consistently verified by the numerical examples presented above..
The generalization of our approach to multidimensional data recovery \cite{ying2024multidimensional,andersson2018esprit} is straightforward.
For future work, we will be investigating ways to optimize sampling locations.

\section*{Declarations}
\textbf{Conflict of interest}\quad The authors declare to have no conflict of interests.

{ 
\bibliographystyle{siam} 
\bibliography{SparseRecovery}
}

\end{document}